\documentclass[12pt]{article}
\usepackage{amssymb,amsfonts}
\def\<{\langle}
\def\>{\rangle}

\def\d{\delta}

\def\btl{\blacktriangleleft}
\def\D{\Delta}

\def\i{\iota}

\def\r{\rho}
\def\lr{\longrightarrow}

\def\o{\otimes}
\def\om{\omega}

\def\vp{\varphi}

\def\s{\sum}
\def\si{\sigma}
\def\v{\varepsilon}

\date{}
\date{}
\begin{document}
\renewcommand{\baselinestretch}{1.2}
\renewcommand{\arraystretch}{1.0}
\title{\bf Quasitriangular ($G$-cograded) multiplier Hopf algebras}
\author {{\bf L. Delvaux$^{a}$, A. Van Daele$^{b}$
 and Shuanhong Wang$^{b, c}$}\\
\small {${}^a$ Department of Mathematics, L. U. C.,
 Universiteitslaan,}\\
\small { B-3590 Diepenbeek, Belgium }\\
\small {E-mail: lydia.delvaux@luc.ac.be}\\
\small {${}^b$ Department of Mathematics, K.U. Leuven,
 Celestijnenlaan 200B, }\\
\small {B-3001 Heverlee, Belgium }\\
\small {E-mail: Alfons.Van Daele@wis.kuleuven.ac.be}\\
\small {${}^b$ Department of mathematics, Southeast University,
 Nanjing 210096, China }\\
\small {E-mail: shuanhwang2002@yahoo.com}}
 \maketitle
\begin{abstract}
 We put the known results on the antipode of a
 usual quasitriangular Hopf algebra into the
 framework of multiplier Hopf algebras. We illustrate with
 examples which can not be obtained by using classical Hopf
 algebras. The focus of the present paper lies on
 the class of the so-called $G$-cograded multiplier Hopf algebras.
 By doing so, we bring the results
 of quasitriangular Hopf group-coalgebras (as introduced by
 Turaev) to the more general framework of multiplier Hopf
 algebras.
\end{abstract}

\vskip 0.6cm
\section*{Introduction}
\def\theequation{0. \arabic{equation}}
\setcounter{equation} {0} \hskip\parindent

The motivating example for quasitriangular Hopf algebras is given
 by the Hopf algebra $H=U_q(g)$, for $g$ a finite-dimensional
 semisimple Lie algebra over $k=\mathbb{C}$, see [Dr1]. However
 $H$ is not quasitriangular in the strict sense of the definition.
 The $R$-matrix lies in a completion of $H\o H$ rather than in
 $H\o H$ itself. The Hopf algebra $H=U_q(g)$ is called "topologically"
 quasitriangular. We believe that an approach with multiplier Hopf
 algebras can evade this problem in purely algebraic terms. The
 notion of a quasitriangular multiplier Hopf algebra is introduced
 in [Z].
\\

 {\bf Definition} [Z]. A regular multiplier Hopf algebra $(A, \D )$
 is called quasitriangular if there exists an invertible
 multiplier $R$ in $M(A\o A)$ which is subject to

 (1) $(\D \o \i )(R)=R^{13}R^{23} \quad (\i \o \D )(R)=R^{13}R^{12}$

(2) $R\D (a)=\D ^{cop}(a)R$ \quad for all $a \in A$

(3) $(\v \o \i )(R)=(i\o \v )(R)=1\in M(A)$.
\\

 Observe that in the equations (1) and  (3) we extended
 a non-degenerate algebra homomorphism to the multiplier
 algebra. This technique is natural in the
 framework of multiplier Hopf algebras. For more details
 on these extensions, we refer to [VD1-Appendix A4].

The present paper is organized as follows.

In the Preliminaries, we review some recent results on
 multiplier Hopf algebras which are used throughout this
 paper.

 In Section 2 we prove that
 the properties on the antipode
 of a quasitriangular Hopf algebras (as given in [Dri 2]) also
 hold for the antipode of a quasitriangular multiplier Hopf
 algebra. The general result is stated in Proposition 2.6. For a
 quasitriangular discrete multiplier Hopf algebra $A$ we can
 express the inner automorphism $S^4$ by using the modular
 multiplier in $A$ and in the reduced dual $\hat{A}$, see
 Proposition 2.9.  To finish Section 2, we give infinite-dimensional
 examples which can not be obtained in the framework of usual Hopf
 algebras.

 In Section 3 we consider a $G$-cograded multiplier Hopf
 algebra $A$ elaborated with a crossing of $G$ on $A$,
 see Preliminaries. In Definition 3.1, we define when
 a $G$-cograded multiplier Hopf algebra with a crossing
 $\pi $ is "$\pi $-quasitriangular".
 Throughout Section 3 we investigate the antipode of a
 $G$-cograded multiplier Hopf algebra which is $\pi
 $-quasitriangular. By using the crossing of $G$ on $A$
  in an appropriate way, we prove that $S^4$ is again
  an inner automorphism on $A$, induced by a grouplike
  multiplier in $M(A)$. The proofs make use of natural
  tools and techniques of multiplier Hopf algebra theory.
  An important class of $\pi $-quasitriangular multiplier Hopf algebras
  is given by the Drinfel'd double constructions for $G$-cograded
  multiplier Hopf algebras, see Theorem 3.11.1.
  We notice that the results of Section 3 apply to the so-called
   Hopf group-coalgebras as introduced by Turaev in [T]. Indeed,
   in [A-De-VD, Theorem 1.5], we have explained how to consider
   a Hopf group-coalgebra in the framework of group-cograded
   multiplier Hopf algebras. Therefore, the result
   of Zunino in  [Zun, Theorem 5.6] is a special case of
   Theorem 3.11.1.

\vskip 0.6cm
\section*{1. Preliminaries}
\def\theequation{0. \arabic{equation}}
\setcounter{equation} {0} \hskip\parindent

Several important constructions and conceptions for classical Hopf
 algebras have their best examples in the finite-dimensional
 cases, e.g. dualizing, Drinfel'd double, quasitriangularity,
 $\cdots$.
 In [VD1], A. Van Daele introduced the much larger class of multiplier
  Hopf algebras. In this framework,
  one considers an algebra $A$ over the field $\mathbb{C}$, with
  or without an identity, but with a non-degenerate
  multiplication map $m$ (viewed as bilinear form)
  on $A\times A$. There is a homomorphism $\D $ from $A$ to the
  multiplier algebra $M(A\o A)$ of $A\o A$. Certain conditions
  on $\D $ (such as coassociativity) are imposed. The motivating
  example in [VD1] is the case where $A=K(G)$, the algebra of
  complex valued, finitely supported functions on an arbitrary
  group $G$. In $A=K(G)$, one uses the pointwise multiplication
  while the coproduct is given by the formula
   $(\D f)(s, t)=f(st)$ with $s, t\in G$ and $f\in A$.

   If $A$ has an identity, one has a usual Hopf algebra. Also for these
 multiplier Hopf algebras, there is a natural notion of left
 and right invariance for linear functionals (called integrals in Hopf
 algebra theory). We only consider regular multiplier Hopf
 algebras (i.e. with a bijective antipode). If $(A, \D )$ has
 invariant functionals, one can construct, in a canonical way, the
 restricted dual $(\hat{A}, \hat{\D} )$.
 We notice that $(\widehat{A}, \widehat{\D} )$ is in the
 same category and $\widehat{\widehat{A}}$ is isomorphic to $A$.
 The duality theorems are investigated in [VD2].
 It is shown that these theorems coincide
 with the usual duality for finite-dimensional Hopf algebras.
 In fact, the evaluation $\< \hat{A}, A\>$ is a special case of a pairing
 of two multiplier Hopf algebras.

  The Drinfel'd double
 construction, which is in the Hopf algebra framework only
 considered for finite-dimensional cases, can be
 done for any pairing $\< A, B\>$ of two multiplier Hopf
 algebras.  Essentially, two multiplier Hopf algebras $A$ and $B$ are paired
 if the product of $A$ (resp. coproduct of $A$) is dual to the
 coproduct of $B$ (resp. product of $B$). Regularity conditions
 are needed to write down this idea in a correct way.
 For details, we refer to [Dra-VD], [De-VD1].
 The Drinfel'd double $D=A\bowtie B^{cop}$ of a pair $\< A, B\>$
 is a twisted tensor product algebra, installed on the tensor product $A\o
 B$. It is proven in [Dra-VD] that the following twist map $T$
 is well-defined. The map $T: B\o A\lr A\o B$ is given by the
 formula $T(b\o a)=\s \< a_{(1)}, S^{-1}(b_{(3)})\>\< a_{(3)}, b_{(1)}\>
 a_{(2)}\o b_{(2)}$ for all $a\in A$ and $b\in B$. The
 comultiplication in $D=A\bowtie B^{cop}$ is given as follows.
 Take $a\in A$ and $b\in B$ $\D (a\bowtie b)=\D (a)\D ^{cop}(b)$
 considered as multiplier in $M(D\o D)$.
 Of course we have $\D ^{cop}(a\bowtie b)=\D ^{cop}(a)\D (b)$.
\\
\\

 Recently, we considered the so-called $G$-cograded multiplier
 Hopf algebras, see [A-De-VD]. This class of
 multiplier Hopf algebras is a generalization of the Hopf
 group-coalgebras as introduced by Turaev in [T]. Let $G$ be any group.
\\

{\bf Definition }[A-De-VD].  A multiplier
 Hopf algebra $B$ is called $G$-cograded if we have

(1) $B=\s _{p\in G}\bigoplus B_p$ with $(B_{p})_{p\in G}$ a family
 of subalgebras of $A$ such that $B_pB_q=0$ if $p\neq q$,

 (2) $\D (B_{pq})(1\o B_q)=B_p\o B_q$ and $(B_p\o 1)\D (B_{pq})
 =B_p\o B_q$ for all $p, q\in G$.
\\
\\

Let $Aut(B)$ denote the group of algebra automorphism on $B$.

{\bf Definition} [De-VD3]. By an action of the group $G$ on $B$,
 we mean a group homomorphism
 $\pi : G\lr Aut(B)$. We require that for all $p\in G$

(1) $\pi _p$ respects the comultiplication on $B$ in the sense
 that for all $b\in B$ $\D (\pi _p(b))=(\pi _p\o \pi _p)(\D (b))$.
\\We call this action {\it admissible} if there is an action $\r $ of $G$
 on itself so that

(2) $\pi _p(B_q)=B_{\r _p(q)}$

(3) $\pi _{\r _p(q)}=\pi _{pqp^{-1}}$
\\for all $p, q$ in $G$.

An admissible action is called a {\it crossing} if furthermore
 $\r _p(q)=pqp^{-1}$ for all $p, q\in G$.
\\
\\

{\bf Theorem} [De-VD3]. Let $B$ be a $G$-cograded multiplier Hopf
 algebra and let $\pi $
 be an admissible action of $G$ on $B$. Then we have a deformation
 $\widetilde{B}$ of $B$ in the following way. As algebra,
 $\widetilde{B}=B$. The comultiplication $\widetilde{\D }$ on
 $\widetilde{B}$ is defined as follows

$$
\widetilde{\D }(b)(1\o b')=(\pi _{q^{-1}}\o \i )(\D (b)(1\o b'))
$$
where $b\in B$ and $b'\in B_q$. The counit $\widetilde{\v }$ is
 the original counit $\v $. The
 antipode $\widetilde{S}$ is given by the formula $\widetilde{S}(b)=
 \pi _{p^{-1}}(S(b))$ for $b\in B_p$.
\\
\\

 There is a Drinfel'd double construction in the framework
 of $G$-cograded multiplier Hopf algebras in the following sense.

{\bf Proposition} [De-VD3].  Let $\< A, B\>$ be a pairing of two
 regular multiplier Hopf algebras. Suppose that $B$ is $G$-cograded.
 Then there exists subspaces $\{A_p\}_{p\in G}$
 of $A$ such that

 (1) $A=\bigoplus _{p\in G}A_p$ and $A_pA_q\subseteq A_{pq}$,

 (2) $\<A_p, B_q\>=0$ whenever $p\neq q$,

 (3) $\<\D (A_p), B_q\o B_r\>=0$ if $q\neq p$ or $r\neq p$
 where $p, q, r\in G$.
\\
\\

{\bf Theorem} [De-VD3]. Let $\< A, B\>$ be a pair of multiplier
 Hopf algebras and assume that $B$ is a $G$-cograded multiplier
 Hopf algebra. Let $\pi $ be a admissible action of $G$ on $B$. The space
 $D^{\pi }=A^{cop}\bowtie \widetilde{B}$ is a multiplier Hopf
 algebra, called the Drinfel'd double, with the multiplication,
 the comultiplication, the counit and the antipode, depending
 on the pairing as well as on the action $\pi $ in the following
 way
\\

$\bullet $ $(a\bowtie b)(a'\bowtie b')=(m_A\o m_B)(\i \o T\o \i
 )(a\o b\o a'\o b')$ where $T(b\o a')=\s \< a'_{(1)}, S^{-1}(b_{(3)})\>
 \< a'_{(3)}, \pi _{p^{-1}}(b_{(1)})\>a'_{(2)}\o b_{(2)}$ for all
 $a'\in A$ and $b\in B_p$.
\\

$\bullet $ $\D (a\bowtie b)=\D ^{cop}(a)\widetilde{\D} (b)$
 where $\D ^{cop}(a)$ and $\widetilde{\D}(b)$ are considered
 as multipliers in $M(D^{\pi }\o D^{\pi })$.
\\

$\bullet $ $\v (a\bowtie b)=\v (a)\v (b)$
\\

$\bullet $ $S(a\bowtie b)=T(\pi _{p^{-1}}(S(b))\o S^{-1}(a))$.

\vskip 0.6cm
\section*{2. Quasitriangular multiplier Hopf algebras}
\def\theequation{0. \arabic{equation}}
\setcounter{equation} {0} \hskip\parindent

Let $A$ be any regular multiplier Hopf algebra. Let $R$
 denote a multiplier in $M(A\o A)$ so that for all $a\in A$ we
 have $(a\o 1)R$ and $R(a\o 1)$ are elements in $A\o M(A)$. We use
 the following notation

$$
(a\o 1)R=\s aR^{(1)}\o R^{(2)}\quad \mbox {and} \quad R(a\o 1)=\s
 R^{(1)}a\o R^{(2)}
$$
in $A\o M(A)$.

To a multiplier $R$ in $M(A\o A)$ which satisfies the above
 conditions, we associate a left multiplier in the sense of the
 following definition.
\\
\\

{\bf 2.1 Definition.} Take $R$ as above. We define
 the left multiplier $u$ as follows. Take $a\in A$. We set

$$
ua=\s S(R^{(2)})R^{(1)}a \quad \, \, \, \, \mbox {in}\, \, A.
$$

It is easy to see that $u$ is a left multiplier of $A$.

Notice that we use the extension to $M(A)$ of the anti-isomorphism
 $S$ on $A$.
\\
\\

In the following lemma we give a condition on $R$ in order that
 $u$ is a two-sided multiplier in $M(A)$.

{\bf 2.2 Lemma.} Take $R\in M(A\o A)$ as above. Suppose that for
 all $a\in A$, $R\D (a)=\D ^{cop}(a)R$ in $M(A\o A)$. Then we have
\smallskip

 (1) $u$ is a multiplier in $M(A)$ so that
 $ua=S^2(a)u$ for all $a\in A$.
\smallskip

 (2) If $R$ is invertible in $M(A\o A)$, then $u$ is invertible in
 $M(A)$. Now $S^2$ is an inner automorphism. Take $a\in A$,
 we have $ S^2(a)=uau^{-1}=S(u)^{-1}a S(u)$.
\\

{\bf Proof.} (1) The original proof of [Dri2] can be easily
 modified into the framework of multiplier Hopf algebras.
 Take any $a, b, c$ in $A$. We have in $A\o A\o A$

$$
\s R^{(1)}a_{(1)}b\o R^{(2)}a_{(2)}\o ca_{(3)} =\s
a_{(2)}R^{(1)}b\o a_{(1)}R^{(2)}\o ca_{(3)}.
$$

Therefore,

$$
\s S^2(ca_{(3)})S(R^{(2)}a_{(2)})R^{(1)}a_{(1)}b=\s
 S^2(ca_{(3)})S(a_{(1)}R^{(2)})a_{(2)}R^{(1)}b.
$$

So,

$$
\s S^2(c)S(R^{(2)})R^{(1)}ab=\s
 S^2(c)S^2(a)S(R^{(2)})R^{(1)}b.
$$

This means that $u(ab)=S^2(a)(ub)$ for all $a, b\in A$. So
 $u$ is a multiplier in $M(A)$ and the right multiplication
 with $u$ is given by the formula $S^2(a)u=ua$ for all $a\in A$.
\\

(2) Suppose that $R$ is invertible in $M(A\o A)$. Let
 $V\in M(A\o A)$ so that $RV=VR=1\o 1$ we suppose again
 that for all $a\in A$, $V(a\o 1), (a\o 1)V$ are in $A\o M(A)$.
 Define the left multiplier $t$ of $A$ as follows
\\

Take $a\in A$, we set $ ta=\s S^{-1}(V^{(2)})V^{(1)}a$.
\\

We can prove in a similar way as above that $t$ is a multiplier in
 $M(A)$. More precisely, we have for all $a\in A$ that
 $ta=S^{-2}(a)t$.

We now calculate that $ut=1$ in $M(A)$. It is sufficient
 to prove that these multipliers equal as left multipliers
 on $A$. Take  any $a\in A$, we have

\begin{eqnarray*}
(ut)a&=&u(\s S^{-1}(V^{(2)})V^{(1)}a)=\s S(V^{(2)})uV^{(1)}a\\
&=&\s S(V^{(2)})S(R^{(2)})R^{(1)}V^{(1)}a=a.
\end{eqnarray*}

Furthermore, $ut=S^2(t)u=1$ and we conclude that
 $u^{-1}=t=S^2(t)$.

From (1), we easily deduce that $S^2$ is inner and for all
 $a\in A$, $S^2(a)=uau^{-1}$. By applying the antipode on both
 sides of this equation and using that $S(A)=A$, we easily
 become the second expression for $S^2$. $\blacksquare$
\\
\\

We now suppose that $A$ is quasitriangular in the sense of [Z].
 We use [Z, Definition 1] in a slightly different formulation.
\\
\\

{\bf 2.3 Definition.} A regular multiplier Hopf algebra is
 called quasitriangular if there is an invertible multiplier
 $R\in M(A\o A)$ so that

 (1) $R\D (a)=\D ^{cop}(a)R$  for all  $a \in A$,

 (2) $(\D \o \i )(R)=R^{13}R^{23}$

 (3) $(\i \o \D )(R)=R^{13}R^{12}$.

We assume furthermore that for all $a$, $R(a\o 1)$ and
 $(a\o 1)R$ are in $A\o M(A)$.
\\

 The following results are known for usual quasitriangular
 Hopf algebras.
\\
\\

 {\bf 2.4 Lemma.} Let $R$ be a generalized $R$-matrix for $A$ as
 in Definition 2.3. Then we have for all $a\in A$

(1) $(\v \o \i )(R)=1=(\i \o \v )(R)$ in $M(A)$

(2) $R^{-1}(a\o 1)=(\i \o S^{-1})(R(a\o 1))$ \quad \, \, $(a\o
 1)R^{-1}=(\i \o S^{-1})((a\o 1)R)$

\quad \, $R^{-1}(a\o 1)=(S\o \i )((S^{-1}(a)\o 1)R)$  \quad $(a\o
1)R^{-1}=(S\o \i )(R(S^{-1}(a)\o 1))$

(3) $(S\o S)(R)=R$ in $M(A\o A)$.
\\
\\

{\bf Proof.} (1) This proof is the same as for usual Hopf
algebras.

(2) We prove the expressions for $R^{-1}$ as left multiplier.
 The other proofs are similar. Take $a, b\in A$. Then
 we have

\begin{eqnarray*}
&&R((\i \o S^{-1})(R(a\o 1)))\\
&=&(\i \o m)(\i \o \i \o S^{-1})(((\i \o \D ^{cop})(R))(a\o 1\o 1))\\
&=&a\o 1.
\end{eqnarray*}

Therefore, $R^{-1}(a\o 1)=(\i \o S^{-1})(R(a\o 1))$
 because $R$ is supposed to be invertible.

Similarly, we calculate

\begin{eqnarray*}
&&(b\o 1)R((S\o \i )((S^{-1}(a)\o 1)R))\\
&=&(m\o \i )(\i \o S\o \i )((b\o S^{-1}(a)\o 1)
((\D \o \i )(R)))\\
&=& ba\o 1.
\end{eqnarray*}

Therefore, $R^{-1}(a\o 1)=(S\o \i )((S^{-1}(a)\o 1)R)$.

\begin{eqnarray*}
\mbox (3) \quad \quad &&((S\o S)(R))(S(a)\o 1)=(S\o S)((a\o
1)R)\quad
 \quad \quad \\
&=&(\i \o S)((S\o \i )((a\o 1)R))=(\i \o S)(R^{-1}(S(a)\o 1))\\
&=&(\i \o S)((\i \o S^{-1})(R(S(a)\o 1)))=R(S(a)\o 1).
\end{eqnarray*}

In the case of a usual Hopf algebra, we get the known results
 when we put $a=1.\quad \blacksquare$
\\
\\

{\bf 2.5 Corollary.} Let $R$ be a generalized $R$-matrix for
 $A$, in the sense of Definition 2.3. Let $u$ be in $M(A)$,
 defined by the formula $ua=\s S(R^{(2)})R^{(1)}a$ for all
 $a\in A$. Then $u$ is invertible in $M(A)$ and $u^{-1}$
 is given as
\begin{eqnarray*}
&&u^{-1}a=\s S^{-2}(R^{2})R^{(1)}a=\s S^{-1}(R^{2})S(R^{(1)})a\\
&&\quad \quad =\s R^{(2)}S^2(R^{(1)})a.
\end{eqnarray*}
\\

{\bf Proof.} Combine the proof of Lemma 2.2(2) and Lemma 2.4 (2).
 $\blacksquare$
\\
\\

Let $A$ be  a quasitriangular multiplier Hopf
 algebra. For all $a\in A$, we have
 $S^2(a)=uau^{-1}=S(u)^{-1}aS(u)$, see Lemma 2.2.
Therefore, we also have
 $$
S^4(a)=uS(u)^{-1}aS(u)u^{-1}, \quad \mbox {for all} \quad a\in A.
 $$

We now prove that $uS(u)^{-1}$ is a
 grouplike multiplier in $M(A)$. Our proof is slightly different
 from the original proof in [Dri2]. We only make use of the
 multiplication by multipliers (instead of using actions),
  combined by the Yang-Baxter equation
 $$
R^{12}R^{13}R^{23}=R^{23}R^{13}R^{12} \, \, \mbox {in}\,  \,
 M(A\o A\o A).
 $$
For a proof of the Yang-Baxter equation for
 quasitriangular multiplier Hopf algebras, we refer to [Z, Proposition
 3].
\\

{\bf 2.6 Proposition.} Let $(A, R)$ be a quasitriangular
 multiplier Hopf algebra. Then for all $a\in A$

 $$
S^4(a)=gag^{-1}
 $$

where $g=uS(u)^{-1}$ is a grouplike multiplier in $M(A)$.
\\

{\bf Proof.} Let $\si $ denote the flip map on
 $A\o A$ and extend it to $M(A\o A)$.
 We calculate the multiplier $\D (u)\si (R)R$.
 Take any $a, b\in A$. Take $a_i, b_i\in A$ so that
 $a\o b=\s \D (a_i)(1\o b_i)$. Then we have
\begin{eqnarray*}
&&(\D (u)\si (R)R)(a\o b)=\s (\D (u)\si (R)R)(\D (a_i)(1\o b_i))\\
&=&\s \D (ua_i)\si (R)R(1\o b_i)=\s \D (S(R^{(2)})R^{(1)}a_i)\si
 (R)R(1\o b_i)\\
&=&\s \D (S(R^{(2)}))\D (R^{(1)}a_i)\si (R)R(1\o b_i)\\
&=&\s \D (S(R^{(2)}))\si (R)R\D (R^{(1)}a_i)(1\o b_i)\\
&=&\s ((\si \circ (S\o S))\D (R^{(2)}))\si (R)R\D (R^{(1)}a_i)(1\o
b_i)
\end{eqnarray*}
Throughout the proof, the letters $R, T, U, V, W, Z$ are used to
 denote the $R$-matrix in $M(A\o A)$.
\begin{eqnarray*}
&=&\s ((\si \circ (S\o S))(R^{(2)}\o U^{(2)}))\si (R)R\D
 (U^{(1)}R^{(1)}a_i)(1\o b_i)\\
&=&\s ((\si \circ (S\o S))(R^{(2)}W^{(2)}\o U^{(2)}V^{(2)}))\si
(R)R(U^{(1)}R^{(1)}a\o V^{(1)}W^{(1)}b)\\
&=&\s (S(\underline
 {U^{(2)}V^{(2)}})T^{(2)}\underline {Z^{(1)}U^{(1)}}R^{(1)}a)\o
 (S(R^{(2)}W^{(2)})T^{(1)}\underline {Z^{(2)}V^{(1)}}W^{(1)}b).
\end{eqnarray*}
If we apply the Yang-Baxter equation to the underlined expressions
 we get
\\
\\
$\s (S(V^{(2)}U^{(2)}) T^{(2)}U^{(1)}Z^{(1)}R^{(1)}a)\o
 (S(R^{(2)}W^{(2)})T^{(1)}V^{(1)}Z^{(2)}W^{(1)}b)$
\\

We notice that for all $a\in A$, we have
\begin{eqnarray*}
&&\s U^{(1)}R^{(1)}a\o S(R^{(2)})U^{(2)}= U^{(1)}R^{(1)}a\o
 S(S^{-1}(U^{(2)})R^{(2)})\\
&=&(\i \o S)(\s U^{(1)}R^{(1)}a\o S^{-1}(U^{(2)})R^{(2)})=a\o 1.
\end{eqnarray*}

Now the above expression can be simplified as follows.
\\
\\
$\s S(U^{(2)})U^{(1)}a\o S(W^{(2)})W^{(1)}b=(u\o u)(a\o b)$.
\\

We obtain $\si (R)R\D (u)=\D (u)\si (R)R=u\o u$ in $M(A\o A)$.\\

Therefore,

\begin{eqnarray*}
&&\D (S(u))=(\si \circ (S\o S))\D (u)\\
&=&(\si \circ (S\o S))((\si (R)R)^{-1}(u\o u))\\
&=&(S(u)\o S(u))\si (((S\o S)(R)(S\o S)(\si (R)))^{-1})\\
&=&(S(u)\o S(u))\si ((R\si (R))^{-1})\\
&=&(S(u)\o S(u))(\si (R)R)^{-1}.
\end{eqnarray*}

Now, it is easy to see that $uS(u)^{-1}$ is grouplike.
 $\blacksquare $
\\
\\
\\

To finish this section, we  suppose that $A$ is a discrete
 multiplier Hopf algebra, as introduced in [VD-Z]. A multiplier
 Hopf algebra $A$ is discrete if $A$ contains
 cointegrals. A left cointegral in $A$ is a non-zero element
 $h\in A$ such that $ah=\v (a)h$ for all $a\in A$. Right
 cointegrals are defined in a similar way.

 In [VD-Z, Theorem 2.10] is proven that a discrete multiplier
 Hopf algebra $A$ always contains non-zero left (and right)
 invariant functionals. Recall that a non-zero linear map $\vp $ on $A$
 is called left invariant if $(\i \o \vp )\D (a)=\vp (a)1$ for all
 $a\in A$. Similarly, non-zero right invariant functionals are
 considered. The reduced dual $\widehat{A}$ is given as the
 following subspace of the full linear dual.

$$
\widehat{A}=\{\vp (a\cdot )\mid a\in A\}=\{\vp (\cdot a)\mid a\in
 A \}.
$$

For more details on this duality, we refer to [VD2]. In the
 present paper, we make use of the so-called modular multiplier
 $\d _A$ which is assigned to a multiplier Hopf algebra $A$ with
 non-zero left (right) invariant functionals. We have from
 [VD2, Propositions 3.8, 3.9] that the modular multiplier
 $\d _A$ in $M(A)$ is uniquely determined such that
 for all $a\in A$

\begin{eqnarray*}
&& (\vp \o \i )\D (a)=\vp (a)\d _A \quad \quad \quad \d _A \,
 \mbox {is invertible and }\\
&& \D (\d _A)=\d _A\o \d _A \quad \quad \v (\d _A)=1\, \quad
 \quad  S(\d _A)=\d ^{-1}_A.\\
&&
\end{eqnarray*}

The evaluation map $\< \widehat{A}, A\>$ is an important case of a
 pairing between multiplier Hopf algebras. According to this pairing,
 we consider the right action of $\widehat{A}$ on $A$, denoted as
 $A\btl \widehat{A}$ which is defined in the following way.
 Take $a\in A$ and $g\in \widehat{A}$, then
 $a\btl g=\s \< g, a_{(1)}\>a_{(2)}$. Observe that
 $a_{(1)}$ is covered by $g$ through the pairing.
 The action $A\btl \widehat{A}$ is unital and extends
 to the multiplier algebra in an obvious way.
\\
\\

The Lemmas 2.7 and 2.8 give general results on discrete multiplier
 Hopf algebras which are used in Proposition 2.9.

 {\bf 2.7 Lemma.} Let $h$ be a left cointegral in $A$. Then
 the evaluation $\widehat{\vp }_h=\< \cdot , h\>$ is a left
 integral on $\widehat{A}$. Furthermore, $ha=\< \d _{\widehat{A}}, a\>h$
 for all $a\in A$.
\\

{\bf Proof.} Clearly $\widehat{\vp} _h$ is a linear functional on
 $\widehat{A}$. For all $f\in \widehat{A}$, we calculate the
 multiplier $(\i \o \widehat{\vp} _h)\D (f)$ in $M(\widehat{A})$.
 Take $g\in \widehat{A}$, then by definition,

 $$
((\i \o \widehat{\vp} _h)\D (f))g=(\i \o \widehat{\vp} _h)(\D
(f)(g\o 1))=\s f_{(2)}(h)f_{(1)}g.
 $$
For $a\in A$, we have

\begin{eqnarray*}
&&(\s f_{(2)}(h)f_{(1)}g)(a)=\s
f_{(2)}(h)f_{(1)}(a_{(1)})g(a_{(2)})\\
&=&\s f(a_{(1)}h)g(a_{(2)})=f(h)g(a).
\end{eqnarray*}
\\
Therefore, $((\i \o \widehat{\vp} _h)\D (f))g=\widehat{\vp}
 _h(f)g$ and this means that $\widehat{\vp} _h$ is a left integral
 on $\widehat{A}$.

 As $\widehat{\vp} _h$ is a left integral on $\widehat{A}$, we
 have that for all $f\in \widehat{A}, (\widehat{\vp} _h\o \i )\D (f)
 =\widehat{\vp} _h(f)\d _{\widehat{A}}$ where $\d _{\widehat{A}}$
 is the modular multiplier in $\widehat{A}$ as defined in [VD2].
 Take any $g\in \widehat{A}$, then we have $g=\vp (x\cdot )$ where
 $\vp $ is a left integral on $A$ and $x\in A$. The equation
 $$
\s \< f_{(1)}, h\>f_{(2)}g=\< f, h\>\d _{\widehat{A}}g
 $$
yields for all $f, g\in \widehat{A}$. By using $(1\o A)\D (A)=A\o
 A$, we easily become that $ha=\< \d _{\widehat{A}}, a\>h$ for all
 $a\in A$. $\blacksquare $
\\
\\

{\bf 2.8 Lemma.} Let $A$ be any discrete multiplier Hopf algebra.
 Let $h$ denote a left cointegral in $A$ and $\d _A$ (resp. $\d _{\widehat{A}}$)
 is the modular multiplier in $M(A)$ (resp. $M(\widehat{A})$). For
 all $a\in A$, we have

(1) $(1\o a)\D (h)=(S(a)\o 1)\D (h)$

(2) $\D (h)(a\o 1)=\D (h)(1\o S(a\btl \d _{\widehat{A}}))$

(3) $\D (h)(a\o 1)=\< \d _{\widehat{A}}, \d _{A}\>^{-1}\s
 S^2(h_{(2)})a\o h_{(1)}\d _A$.

Notice that the modular $\< \d _{\widehat{A}}, \d _A\>$ has a
 meaning because $\d _A$ is a grouplike multiplier.
\\

 {\bf Proof.} (1) See [VD-Z, Proposition 2.5].

(2) Take any $a, b$ in $A$. Then we have

\begin{eqnarray*}
\D (h)(a\o b)&=&\s h_{(1)}a_{(1)}\o h_{(2)}\v (a_{(2)})b\\
&=&\s h_{(1)}a_{(1)}\o h_{(2)}a_{(2)}S(a_{(3)})b
=\s \D (ha_{(1)})(1\o S(a_{(2)})b)\\
&=&\s \< \d _{\widehat{A}}, a_{(1)}\>
 \D (h)(1\o S(a_{(2)})b)=\D (h)(1\o S(a\btl \d _{\widehat{A}})b).
\end{eqnarray*}

Now the statement follows.

(3) Let $\psi $ be a right integral on $A$ so that
 $\psi (h)=1$. Apply the operator $\psi \o \i $ on the both sides
 of the equation in (1). Then we get for all $a\in A$

$$
a=\s \psi (S(a)h_{(1)})h_{(2)}
$$
and equivalently,

$$
S(a)=\s \psi (ah_{(1)})S^2(h_{(2)}).
$$

We put $\psi =\vp \circ S$ for some left integral on $A$. Then
 we have $\vp (h)=\< \d _{\widehat{A}}, \d _A\>^{-1}$ because
 $\psi (a)=\vp (a\d _A)$ for all $a\in A$. Apply the operator
 $\i \o \vp $ to both sides of the equation in (1). Then we get
 for all $a\in A$
$$
\s \vp (ah_{(2)})h_{(1)}=\< \d _{\widehat{A}}, \d _A\>^{-1}S(a).
$$

If we combine these equations, we obtain that
$$
\< \d _{\widehat{A}}, \d _A\>\s \vp (ah_{(2)})h_{(1)} =\s \vp
(ah_{(1)}\d _A)S^2(h_{(2)}).
$$

As $\vp $ is faithful in the sense that $\vp (ab)=0$ for all
 $a\in A$ implies that $b=0$, we have obtain that
$$
\< \d _{\widehat{A}}, \d _A\>\s h_{(2)}\o h_{(1)}b =\s h_{(1)}\d
 _A\o S^2(h_{(2)})b
$$
for all $b\in A$. Now the statement follows. $\blacksquare $
\\
\\

We are ready to calculate the grouplike multiplier $g=uS(u)^{-1}$
 for a discrete multiplier Hopf algebra which is quasitriangular.

 {\bf 2.9 Proposition.} Let $A$ be a discrete quasitriangular
 multiplier Hopf algebra. Let $\d _A$ (resp. $\d _{\widehat{A}}$)
 denote the modular multiplier in $M(A)$ (resp. $M(\widehat{A})$).
 Then we have
 $$
uS(u)^{-1}=\d ^{-1}_A((\i \o \< \d ^{-1}_{\widehat{A}}, \cdot
 \>)(R)) \quad \mbox {in}\, \, M(A).
 $$
\\

{\bf Proof.} In the right hand side of the statement the extension
 of the algebra homomorphism $(\i \o \<\d _{\widehat{A}}^{-1},
  \cdot \>)$ to
 the multiplier algebra $M(A\o A)$ is considered. Throughout the
 proof, $h$ denotes a left cointegral in $A$ and $u$ is the
 multiplier in $M(A)$, associated to $R$-matrix $R$. As $A$ is
 almost co-commutative, we have
 $$
R\D (h)=\D ^{cop}(h)R \quad \mbox {in}\, \, M(A\o A)
 $$
This means that for all $a\in A$ we have
$$
\s R^{(1)}h_{(1)}a\o R^{(2)}h_{(2)}\stackrel {(*)}{=} \s
 h_{(2)}R^{(1)}a\o h_{(1)}R^{(2)}.
$$

By using Lemma 2.8, the left hand side of
 the equation (*) equals

\begin{eqnarray*}
 \s R^{(1)}S(R^{(2)})h_{(1)}a\o h_{(2)} &=&\< \d
 _{\widehat{A}}, \d
 _{A}\>^{-1}S(u)\s S^2(h_{(2)})a\o h_{(1)}\d _A \\
 &=&S(u)\s S^2(h_{(2)})\d _A^{-1}a\o h_{(1)}.
\end{eqnarray*}

By using Lemma 2.8(2), the right hand side of the equation (*) is
 given as

$$
\s h_{(2)}S(R^{(2)}\btl \d _{\widehat{A}})R^{(1)}a\o h_{(1)}.
$$

As $(\i \o \D )(R)=R^{13}R^{12}$, this last expression can be
 written as

$$
\s \< \d _{\widehat{A}}, R^{(2)}\>h_{(2)}uR^{(1)}a\o h_{(1)}.
$$

The equation (*) can be written as

$$
S(u)\s S^2(h_{(2)})\d ^{-1}_Aa\o h_{(1)}=\s \< \d _{\widehat{A}},
 R^{(2)}\>h_{(2)}uR^{(1)}a\o h_{(1)}.
$$

Let $\psi $ be a right integral on $A$ so that $\psi (h)=1$. Apply
 the operator $(\i \o \psi )$ on both sides of the above equation.
 We become that

$$
S(u)\d ^{-1}_Aa=\s \< \d _{\widehat{A}},
 R^{(2)}\>uR^{(1)}a.
$$

Therefore, we have in $M(A)$
$$
S(u)\d ^{-1}_A=u(\i \o \< \d _{\widehat{A}},
 \cdot \>)(R)
$$
and so,
$$
u^{-1}S(u)=((\i \o \< \d _{\widehat{A}},
 \cdot \>)(R))\d _A.
$$

We easily obtain that $uS(u)^{-1}=\d ^{-1}_A(\i \o \< \d
 _{\widehat{A}}^{-1}, \cdot \>)(R)$. $\blacksquare$
\\
\\

{\bf 2.10 Remark.} Take $A$ a finite-dimensional Hopf algebra.
 In [Dri2] and [Ra] is given a similar result as in Proposition
 2.9.
 \\
 \\
 \\

 We now give examples which can not be obtained in the framework
  of usual Hopf algebras.

\vskip 0.3cm

{\it  The dual multiplier Hopf algebra of a Ore-extension}

\vskip 0.3cm

General Ore-extensions of cyclic groups are introduced and
 elaborated in [B-D-G-M]. We consider the special case $A=(kC)_{2, -1,
 1}$.

As algebra, $A$ is generated by an invertible element $a$ and an
 element $b$, so that $b^2=0$. The multiplication in $A$ is
 induced by the commutation rule $ab=-ba$. The coalgebra structure
 is given as follows

\begin{eqnarray*}
&& \D (a)= a\o a, \quad \, \v (a)=1, \quad \, S(a)=a^{-1} \\
&& \D (b)= (a\o b) + (b\o 1), \quad \, \v (b)=0, \,  \quad
 S(b)=-a^{-1}b
\end{eqnarray*}

In [B-D-G-M] is proven that $A$ is co-Frobenius. A left integral
 is given by the functional $\vp $ when $\vp (a^{-1}b)=1$
 and $\vp $ is zero elsewhere on the basis.
 The modular element in $A$ is given as $\d _A=a^{-1}$.
 This means that $(\vp \o \i )\D (a^mb^n)
 =\vp (a^mb^n)a^{-1}$ for all $m\in \mathbb{Z}$ and $n=0, 1$.

 The dual multiplier Hopf algebra $\widehat{A}$, defined as $\widehat{A}=
 \{ \vp (x\cdot )\mid x\in A\}$, is a discrete multiplier Hopf
 algebra, see [VD2]. This means that there are cointegrals in $\widehat{A}$.
 We briefly describe the multiplier Hopf algebra $\widehat{A}$.
 For more details, we refer to [De]. A linear basis of $\widehat{A}$
 is given by the functionals $\om _{m, n}, m\in \mathbb{Z}$ and $n=0, 1$
 where $\om _{m, n}(a^mb^n)=1$ and $\om _{m, n}$ is zero elsewhere
 on the basis of $A$. The product in $\widehat{A}$ is induced by
 the following commutation rules

\begin{eqnarray*}
&& \om _{p, 0}\om _{q, 0}=\d _{p, q}\om _{p, 0} \quad
 \quad
 \om _{p, 0}\om _{q, 1}=\d _{p-q, 1}\om _{q, 1}\\
&& \om _{q, 1}\om _{p, 0}=\d _{q-p, 0}\om _{p, 1} \quad
 \quad \om _{p, 1}\om _{q, 1}=0
\end{eqnarray*}
where $p, q\in \mathbb{Z}$.

The coalgebra structure in $\widehat{A}$ is given as follows.

\begin{eqnarray*}
&& \D (\om _{p, 0})=\s _{q\in \mathbb{Z}}\om _{q, 0}\o \om _{p-q,
 0} \quad \quad  \D (\om _{p, 1})=\s _{r\in \mathbb{Z}, s=0, 1}
 (-1)^{s(p-s)}\om _{r, s}\o \om _{p-r, 1-s}\\
&&\\
&& S(\om _{p, 0})=\om _{-p, 0} \quad \quad S(\om _{p,
 1})=(-1)^p \om _{-p-1, 1}\\
&& \v (\om _{p, 0})=\d _{p, 1}\quad  \quad \v (\om _{p, 1})=0.
\end{eqnarray*}

A left cointegral in $\widehat{A}$ is given by the functional
 $\om _{-1, 1}$, see also Lemma 2.7.

A left integral on $\widehat{A}$ is given by $\widehat{\vp }$
 which is defined on $\widehat{A}$ as follows:

$$
\widehat{\vp }(\om _{p, 0})=0 \, \, \, \mbox {for all} \, \, p\in
\mathbb{Z}, \quad
 \widehat{\vp }(\om _{p, 1})=1 \, \, \, \mbox {for all} \, \,
 p\in \mathbb{Z}.
$$

 The modular multiplier in $M(\widehat{A})$ is given by
  $$
\d _{\widehat{A}}=\s _{p\in \mathbb{Z}}(-1)^p\om _{p, 0}.
  $$

Furthermore, $\widehat{A}$ is quasitriangular. A generalized
 $R$-matrix is given by the multiplier $R$ in $M(\widehat{A}\o \widehat{A})$
 $$
R=\s _{p\in \mathbb{Z}}\d _{\widehat{A}}^p\o \om _{p, 0}.
 $$
One can easily check that

(1) $(\D \o \i )(R)=R^{13}R^{23}$ in $M(\widehat{A}\o
 \widehat{A}\o \widehat{A})$

 (2) $(\i \o \D )(R)=R^{13}R^{12}$ in $M(\widehat{A}\o
 \widehat{A}\o \widehat{A})$

 (3) $R\D (\om _{p, 0})=\D ^{cop}(\om _{p, 0})R$,
 \quad $R\D (\om _{p, 1})=\D ^{cop}(\om _{p, 1})R$.
\\

 It is easy to see that the multiplier $R$ is covered
 by the elements $\om _{k, l}$ ($k\in \mathbb{Z}$ and $l=0, 1$)
 in the sense that $R(\om _{k, l}\o 1)$, $(\om _{k, l}\o 1)R$
 are in $\widehat{A}\o M(\widehat{A})$.
\\

 We have $R^{-1}=R$. The multiplier $u$ in $M(\widehat{A})$
 defining the square of the antipode as an inner
 automorphism (see Lemma 2.2) is given
  by the modular multiplier i.e.,
  $u=\d _{\widehat{A}}$. Therefore, $uS(u)^{-1}=1$ in $M(\widehat{A})$.
  This result indicates that $S^4$ is the identity on
  $\widehat{A}$, see Proposition 2.6. To finish, we notice
  that $\widehat{A}$ is a discrete quasitriangular multiplier
  Hopf algebra.

To illustrate the formula in Proposition 2.9, we calculate in
 $M(\widehat{A})$

\begin{eqnarray*}
\d ^{-1}_{\widehat{A}}((\i \o \<  \d ^{-1}_{A}, \cdot
 \>)(R))&=& (\s _{p\in \mathbb{Z}}(-1)^p\om _{p, 0})
 (\s _{q\in \mathbb{Z}} \d ^{q}_{\widehat{A}}
 \< a, \om _{q, 0}\>)\\
&=& (\s _{p\in \mathbb{Z}}(-1)^p\om _{p, 0})
 \d _{\widehat{A}}=1 \quad \mbox {in} \, \, \, M(\widehat{A}).
\end{eqnarray*}

\vskip 1.5cm {\it  The Drinfel'd double of a Ore-extension}
 \vskip 0.3cm

We start with the Ore-extension $A=(kC)_{2, -1, 1}$. In
 previous example, we considered the dual multiplier Hopf algebra
 $\widehat{A}$. By the construction of $\widehat{A}$, we have that
 $\< \widehat{A}, A\>$ is a pairing of multiplier Hopf algebras in
 the sense of [Dra-VD]. We consider the Drinfel'd  double of
 the pair $\< \widehat{A}, A\>$. Let $D$ denote $D=\widehat{A}\bowtie
 A^{cop}$. For the general construction of the Drinfel'd double
 of a pair of multiplier Hopf algebras, we refer to [Dra-VD] and
 [De-VD1]. One can show that the maps

\begin{eqnarray*}
&& \widehat{A}\lr M(\widehat{A}\bowtie A): \, \, \om _{k,
 l}\mapsto \om _{k, l}\bowtie 1, \quad k\in \mathbb{Z}, \, l=0,
 1,\\
&& A\lr M(\widehat{A}\bowtie A): \, \, a^pb^q \mapsto 1\bowtie
 a^pb^q, \quad p\in \mathbb{Z}, \, q=0, 1
\end{eqnarray*}
are algebra embeddings, so that $\om _{k, l}\bowtie a^pb^q
 =(\om _{k, l}\bowtie 1)(1\bowtie a^pb^q)$ and
 $(1\bowtie a^pb^q)(\om _{k, l}\bowtie 1)=T(a^pb^q\bowtie \om _{k, l})$
 where $T$ is the twist map $T: A\o \widehat{A}\lr \widehat{A}\o
 A$ defined by the formula
$$
T(x\o y)=\s \< y_{(1)}, S^{-1}(x_{(3)})\> \< y_{(3)}, x_{(1)}\>
 y_{(2)}\o x_{(2)}
$$
for all $x\in A$ and $y\in \widehat{A}$.

If we identity $\widehat{A}$ and $A$ with their images in
 $\widehat{A}\o A$, we find that $D=\widehat{A}\bowtie A^{cop}$
 is the algebra generated by $\widehat{A}$ and $A$ with the
 following commutating relations

\begin{eqnarray*}
&& a\om _{p, 0}=\om _{p, 0}a \, \, \, \, \, \quad \quad b\om _{p,
 0}=\om _{p-1, 0}b\\
&& a\om _{p, 1}=-\om _{p, 1}a \quad \quad b\om _{p, 1}= \om
 _{p, 0}-(-1)^p\om _{p, 0}a-\om _{p-1, 1}b
\end{eqnarray*}
where $p\in \mathbb{Z}$.

The comultiplication in $D=\widehat{A}\bowtie A^{cop}$ is given by
 the formula

$$
\D (\om _{p, k}a^rb^l)=\D (\om _{p, k})\D ^{cop}(a^rb^l)
$$
for all $p, r\in \mathbb{Z}$ and $k, l=0, 1$. One can check
 that $R=\s _{p\in \mathbb{Z}, k=0, 1}\om _{p, k}\o a^pb^k$
 in $M(D\o D)$ is a generalized $R$-matrix for $D$. This result
 also follows from [De-VD2, Theorem 4.7]. Remark that for all
 $p\in \mathbb{Z}$ and $k=0, 1$ we have $R(\om _{p, k}\o 1)$
 and $(\om _{p, k}\o 1)R$ in $D\o D$.

 We calculate the associated multiplier $u$ in $M(D)$ in the sense
 of Lemma 2.2.

$$
 u=\s _{p\in \mathbb{Z}, k=0, 1}S^{-1}(a^pb^k)\om _{p, k}
$$

that is,

$$
u=\s _{p\in \mathbb{Z}}a^{-p}\om _{p, 0}+\s _{p\in
 \mathbb{Z}}(-1)^pa^{-1-p}b\om _{p, 1}.
$$

By using the commutation rules in $D$, we obtain

$$
u=\s _{p\in \mathbb{Z}}(-1)^p\om _{p, 0}a^{-p-1}+\s _{p\in
 \mathbb{Z}}\om _{p, 1}a^{-p-2}b.
$$

Following Corollary 2.5, the inverse multiplier $u^{-1}$
 in $M(D)$ is given as

$$
 u^{-1}=\s _{p\in \mathbb{Z}, k=0, 1}a^pb^kS^2(\om _{p, k})
$$

that is,

$$
u^{-1}=\s _{p\in \mathbb{Z}}a^{p}\om _{p, 0}-\s _{p\in
 \mathbb{Z}}a^{p}b\om _{p, 1}.
$$

Again, by using the commutation rules in $D$, we get

$$
u^{-1}=\s _{p\in \mathbb{Z}}(-1)^p\om _{p, 0}a^{p+1}+\s _{p\in
 \mathbb{Z}}(-1)^p\om _{p-1, 1}a^{p}b.
$$

Therefore, we have

$$
S(u)^{-1}=\s _{p\in \mathbb{Z}}(-1)^pa^{-p-1}\om _{-p, 0}
 -\s _{p\in \mathbb{Z}}(-1)^pa^{-1-p}b\om _{-p, 1}.
$$

By the commutation rules in $D$, we obtain

$$
S(u)^{-1}=\s _{p\in \mathbb{Z}}\om _{p, 0}a^p
 -\s _{p\in \mathbb{Z}}\om _{p, 1}a^pb.
$$

The grouplike multiplier $g=uS(u)^{-1}$ which installs
 $S^4$ as an inner automorphism, in the sense of Theorem
 2.6, is given by the expression
\begin{eqnarray*}
g=uS(u)^{-1}&=&\s _{p\in \mathbb{Z}}(-1)^p\om _{p, 0}a^{-1}
 -\s _{p, q\in \mathbb{Z}}(-1)^p\om _{p, 0}a^{-p-1}
 \om _{q, 1}a^qb\\
&&+\s _{p, q\in \mathbb{Z}}\om _{p, 1}a^{-p-2}b
 \om _{q, 0}a^q-\s _{p, q\in \mathbb{Z}}\om _{p, 1}a^{-p-2}b
 \om _{q, 1}a^qb.
\end{eqnarray*}

By using the commutation rules in $D$, we obtain

$$
g=uS(u)^{-1}=\s _{p\in \mathbb{Z}}(-1)^p\om _{p, 0}a^{-1}=\d
 _{\widehat{A}}\bowtie \d _A.
$$
Observe that $g$ is grouplike, but not trivial.

\vskip 0.6cm
\section*{3. Quasitriangular $G$-cograded multiplier Hopf algebras}
\def\theequation{1. \arabic{equation}}
\setcounter{equation} {0} \hskip\parindent

Let $G$ be a group and let $A$ be a $G$-cograded multiplier Hopf
 algebra elaborated with an admissible action $\pi $. We use
 the deformated  multiplier Hopf algebra $(\widetilde{A}, \widetilde{\D})$
 as reviewed in the preliminaries.

{\bf 3.1 Definition.}  Let $A$ be a $G$-cograded multiplier Hopf
 algebra with an admissible action $\pi $. We
 call $A$  $\pi $-quasitriangular
 if there is an invertible multiplier
 $R\in M(A\o A)$ so that

 (1) $(\pi _p\o \pi _p)(R)=R$ for all $p\in G$

 (2) For all $a\in A$, $R\D (a)=\tilde{\D }^{cop}(a)R$
 in $M(A\o A)$

 (3) $(\tilde{\D }\o i)(R)=R^{13}R^{23}$

 (4) $(i\o \D )(R)=R^{13}R^{12}$.

We furthermore suppose that $R(a\o 1), (a\o 1)R\in A\o M(A)$
 for all $a\in A$. We call $R$ the generalized $\pi $-matrix
 for $A$. We denote $R(a\o 1)=\s R^{(1)}a\o R^{(2)}$ in $A\o M(A)$
 and $(a\o 1)R=\s aR^{(1)}\o R^{(2)}$.
\\
\\

If $\pi $ is the trivial action, we recover Definition 2.3. The
 Drinfel'd double construction in the framework of $G$-cograded
 multiplier Hopf algebras is again the prototype of $\pi
 $-quasitriangularity. This statement will be proven in
 Theorem 3.11.1.
 \\
\\

{\bf 3.2 Lemma.} Let $(A, R)$ be a $\pi $-quasitriangular
 $G$-cograded multiplier Hopf algebra. Then we have

 (1) $(\v \o i)(R)=(i\o \v )(R)=1$ in
 $M(A)$

(2) $R^{-1}(a\o 1)=(i\o S^{-1})(R(a\o 1))$

\quad \, $R^{-1}(a\o 1)=(\widetilde{S}\o
i)((\widetilde{S}^{-1}(a)\o
 1)R)$

(3)  $(\widetilde{S}\o S)(R)=R$.
\\

{\bf Proof.} The proofs of this lemma are very similar to
 the proofs of Lemma 2.4.

\begin{flushright}
$\blacksquare$
\end{flushright}

The generalized $\pi $-matrix for a $G$-cograded multiplier Hopf
 algebra satisfies the Yang-Baxter equation in the sense of the
 following proposition.
\\

{\bf 3.3 Proposition.}  Let $A$ be a $\pi $-quasitriangular
 multiplier Hopf algebra. The $\pi $-Yang-Baxter Equation
 for the $\pi $-matrix $R$ is given
 as follows. Take $a\in A$ and $b\in A_p$, for any $p\in G$.
 Then we have
\begin{eqnarray*}
(R^{23}R^{13}R^{12})(a\o b\o 1) &=&(R^{12}(i\o i\o \pi
 _{p^{-1}})(R^{13})R^{23})(a\o b\o 1)\\
&=&(R^{12}(\pi _p\o i\o i )(R^{13})R^{23})(a\o b\o 1).
\end{eqnarray*}

{\bf Proof.} Clearly we have

$$
R^{23}R^{13}R^{12}=R^{23}((i\o \D )(R)) \stackrel {(*)}{=}((i\o
\widetilde{\D}^{cop})(R))(1\o R)
$$

We here explain $(*)$. For all $a, b\in A$

\begin{eqnarray*}
&&((i\o \widetilde{\D}^{cop})(a\o b))(1\o R)
=(a\o 1\o 1)(1\o \widetilde{\D}^{cop}(b)R)\\
&=&(a\o 1\o 1)(1\o R \D (b))=(1\o R)((i\o \D )(a\o b)).
\end{eqnarray*}

Let $U=R=T$ denote the generalized $\pi $-matrix for $A$.
 The Equation (*) means that for all $a\in A$ and $b\in A_p$,
 $p\in G$, we have

\begin{eqnarray*}
(R^{23}R^{13}R^{12})(a\o b\o 1)
&=&\s ((i\o \widetilde{\D}^{cop})(R))(a\o R^{(1)}b\o R^{(2)})\\
&=&\s U^{(1)}a\o (\widetilde{\D}^{cop}(U^{(2)}))(R^{(1)}b\o R^{(2)})\\
&=&\s U^{(1)}a\o ((i\o \pi _{p^{-1}})(\D ^{cop}(U^{(2)})))
(R^{(1)}b\o R^{(2)})\\
&=&\s U^{(1)}T^{(1)}a\o (U^{(2)}\o \pi
 _{p^{-1}}(T^{(2)}))(R^{(1)}b\o R^{(2)})\\
&=&\s U^{(1)}T^{(1)}a\o U^{(2)}R^{(1)}b\o \pi
 _{p^{-1}}(T^{(2)})R^{(2)}.
\end{eqnarray*}
So we have proven the first equation of the statement.
 The second equation easily follows by using condition (1)
 in Definition 3.1. $\blacksquare$
\\
\\
\\

We now consider a $G$-cograded multiplier Hopf algebra $A$. Let
 $\pi $ be a crossing of $G$ on $A$. This means that for all
  $p, q\in G$, we have $\pi _p(A_q)=A_{pqp^{-1}}$.
 Suppose that $R$ is a generalized $\pi $-matrix for $A$.
 We bring the results on the antipode, as proven in Section 2,
 to this framework.

{\bf 3.4 Definition.} Take the assumptions as above. Define
 the left multiplier $\widetilde{u}$ of $A$ as follows

 $$
\widetilde{u}a=\s \widetilde{S}(R^{(2)})R^{(1)}a \quad \quad \, \,
 \mbox {for all}\,\, a\in A.
$$
\\

{\bf 3.5 Remark.} Take $a\in A_p$, then
 $\widetilde{u}a=\s \pi _p(S(R^{(2)}))R^{(1)}a$.
 \\
 \\

{\bf 3.6 Lemma.} Take the notations as above. Take
 $p\in G$ and $a\in A_p$. Then we have

(1) $\widetilde{u}$ is a two-sided multiplier in $M(A)$, the right
 multiplication with $\widetilde{u}$ is given

 \quad \, as follows

 $$
a\widetilde{u}=\widetilde{u}\pi _{p^{-1}}(S^{-2}(a)).
  $$

\quad \, Furthermore, the multiplier $\widetilde{u}$ is invariant
 for each automorphism $\pi _p$.
 \\

(2) $\widetilde{u}$ is invertible in $M(A)$. The square of the
 antipode is given as follows

$$
S^2(a)=\widetilde{u}(\pi _{p^{-1}}(a))\widetilde{u}^{-1}
 $$

\quad \, or equivalently,

$$
S^2(a)=S(\widetilde{u})^{-1}(\pi _p(a))S(\widetilde{u}).
$$
\\

{\bf Proof.} (1) Take $p\in G$. Let $a, b, c\in A_p$

$$
\s R^{(1)}a_{(1)}b\o R^{(2)}a_{(2)}\o ca_{(3)}=\s
a_{(2)}R^{(1)}b\o \pi
 _{p^{-1}}(a_{(1)})R^{(2)}\o ca_{(3)}.
$$

Therefore,

\begin{eqnarray*}
&&\s R^{(1)}a_{(1)}b\o \pi _p(R^{(2)}a_{(2)})\o \pi _p(ca_{(3)})\\
&&\quad \quad =\s a_{(2)}R^{(1)}b\o a_{(1)}\pi _p(R^{(2)})\o
 \pi _p(ca_{(3)}).
\end{eqnarray*}

This implies that

\begin{eqnarray*}
&&\s \pi _p(S^2(ca_{(3)}))\pi _p(S(R^{(2)}a_{(2)}))R^{(1)}a_{(1)}b \\
&&\quad \quad =\s \pi _p(S^2(ca_{(3)}))\pi
 _p(S(R^{(2)}))S(a_{(1)})a_{(2)}R^{(1)}b.
\end{eqnarray*}

The left hand side of the above equation can be written as

\begin{eqnarray*}
&&\s \pi _p(S(R^{(2)}a_{(2)}S(ca_{(3)})))R^{(1)}a_{(1)}b \\
&&\quad \quad =\s \pi
 _p(S(R^{(2)}S(c)))R^{(1)}ab=\pi
 _p(S^2(c))(\widetilde{u}ab).
\end{eqnarray*}

The right hand side of the above equation can be written as

$$
\s \pi _p(S^2(ca))\pi _p(S(R^{(2)}))R^{(1)}b =\pi _p(S^2(c))\pi
 _p(S^2(a)(\widetilde{u}b).
$$

Therefore, we have for all $a, b\in A_p$

$$
\widetilde{u}(ab)=\pi _p(S^2(a))(\widetilde{u}b).
$$

So,  $\widetilde{u}$ is a two-sided multiplier
 in $M(A)$ where the right multiplication with
 $\widetilde{u}$ is given as follows

$$
a\widetilde{u}=\widetilde{u}\pi _{p^{-1}}(S^{-2}(a)) \quad \mbox
 {for all}\, \, a\in A_p.
$$

To finish we remark that $\widetilde{u}$ is invariant
 for each automorphism $\pi _p$ because we have
 $(\pi _p\o \pi _p)(R)=R$.
\\

 (2) We  now prove that $\widetilde{u}$ is invertible in $M(A)$.
 We define the left multiplier $\widetilde{t}$ as follows

 $$
\widetilde{t}a=\s R^{(2)}S^2(R^{(1)})a\quad \mbox {for all}\, \,
 a\in A.
 $$

First one can prove, in a similar way as in (1), that
 $\widetilde{t}$ is a multiplier in $M(A)$.
\smallskip

Notice that for $a\in A_p$, we obtain
 $a\widetilde{t}=\widetilde{t}\pi _p(S^2(a))$.
\smallskip

We now prove that $\widetilde{t}$ is the inverse of
 $\widetilde{u}$ in $M(A)$. Take $a\in A_p$.
 Let $R$ and $V$ denote the $\pi $-matrix in $M(A\o A)$.

\begin{eqnarray*}
\widetilde{t}\widetilde{u}a&=&\s
 R^{(2)}S^2(R^{(1)})\widetilde{u}a=\s
 R^{(2)}\widetilde{u}\pi _{p^{-1}}(R^{(1)}\pi _p(a))\\
&=&\s R^{(2)}\pi _p(S(V^{(2)}))V^{(1)}
\pi _{p^{-1}}(R^{(1)}\pi _p(a))\\
&=&\s S(S^{-1}(V^{(1)}\pi _{p^{-1}}(R^{(1)}\pi _p(a)))
 \pi _p(V^{(2)})S^{-1}(R^{(2)}))\\
&\stackrel {(*)}{=}& S(S^{-1}(a))=a.
\end{eqnarray*}

We explain (*). From Proposition 3.2(2), we have

\begin{eqnarray*}
&&(\s V^{(1)}\o V^{(2)})(R^{(1)}a\o S^{-1}(R^{(2)}))=a\o 1\\
&\Rightarrow &\s V^{(1)}R^{(1)}a\o V^{(2)}S^{-1}(R^{(2)})=a\o 1\\
&\Rightarrow &\s S^{-1}(V^{(1)}R^{(1)}a)\pi _p(V^{(2)})
 S^{-1}(\pi _p(R^{(2)}))=S^{-1}(a)\\
&\Rightarrow &\s S^{-1}(V^{(1)}\pi _{p^{-1}}(R^{(1)}\pi _p(a)))
 \pi _p(V^{(2)})S^{-1}(R^{(2)})=S^{-1}(a).
\end{eqnarray*}

As $\widetilde{t}\widetilde{u}$ and $1$ equal as left multipliers,
 they also equal as multipliers in $M(A)$.

We also have for all $a\in A_p$

\begin{eqnarray*}
a&=&\widetilde{t}\widetilde{u}a=\pi
 _{p^{-1}}(S^{-2}(\widetilde{u}a))\widetilde{t}\\
&=&\pi _{p^{-1}}(S^{-2}(\pi _p(S^2(a))\widetilde{u}))\widetilde{t}
=aS^{-2}(\widetilde{u})\widetilde{t}.
\end{eqnarray*}

Therefore, $S^{-2}(\widetilde{u})\widetilde{t}=1$ and
 also $\widetilde{u}S^2(\widetilde{t})=1$ as multipliers
 in $M(A)$. So we obtained that
 $\widetilde{u}$ is invertible, more precisely,
 $(\widetilde{u})^{-1}=\widetilde{t}=S^2(\widetilde{t})$.

To finish, we observe that the expressions for $S^2$ as given in
 the statement of this lemma easily follow from part (1).
 $\blacksquare$
\\
\\
\\

 By using Lemma 3.6, we easily see that $S^4$ is a usual inner automorphism.
 Indeed, take $a\in A_p$. Then we have

\begin{eqnarray*}
S^4(a)&=&S^2(S^2(a))=\widetilde{u}(\pi _{p^{-1}}(S^2(a)))
 \widetilde{u}^{-1}\\
&=&\widetilde{u}(\pi _{p^{-1}}(S(\widetilde{u})^{-1}\pi _p(a)
 S(\widetilde{u})))\widetilde{u}^{-1}\\
&=&\widetilde{u}S(\widetilde{u})^{-1}aS(\widetilde{u})\widetilde{u}^{-1}.
\end{eqnarray*}

 In Proposition 3.10, we prove the much stronger result that $\widetilde{u}S(\widetilde{u})^{-1}$
 is a grouplike multiplier in $M(A)$.
 First we need to prove some technical lemmas.
\\
\\

{\bf 3.7 Lemma.} Take $a\in A_{pq}$ and $b\in A_q$, then
 we have

\begin{eqnarray*}
&&((i\o \pi _{p^{-1}})\si (R))R\D (a)(1\o b)= \\
&&\quad \quad ((\pi _{pq^{-1}p^{-1}}\o \pi _{p^{-1}})\D (a)) ((i\o
 \pi _{p^{-1}})(\si (R)))R(1\o b).
\end{eqnarray*}

{\bf Proof.} Observe that this result is trivial done if
 $\pi $ is the trivial action of $G$ on $A$.
 Take $a\in A_{pq}$ and $b\in A_q$. Then we have

\begin{eqnarray*}
&& R\D (a)=\widetilde{\D }^{cop}(a)R\\
&\Rightarrow & R\D (a)(1\o b)=\widetilde{\D }^{cop}(a)R(1\o b)\\
&\Rightarrow & R\D (a)(1\o b)=\s (a_{(2)}\o \pi
 _{p^{-1}}(a_{(1)}))R(1\o b)=((i\o \pi _{p^{-1}})\D ^{cop}(a))R(1\o
 b).
\end{eqnarray*}

Therefore,
\begin{eqnarray*}
&&((i\o \pi _{p^{-1}})\si (R))R\D (a)(1\o b)\\
&=&((i\o \pi _{p^{-1}})\si (R\D (a)))R(1\o b)\\
&=&((i\o \pi _{p^{-1}})\si (\widetilde{\D }^{cop}(a)R))R(1\o b)\\
&=&((i\o \pi _{p^{-1}})(\widetilde{\D }(a))((i\o \pi
 _{p^{-1}})(\si (R)))R(1\o b)\\
&&\quad \mbox {we notice that}\, \,
 \pi _{p^{-1}}(a_{(2)})\in A_q\Rightarrow a_{(2)}\in
 A_{pqp^{-1}}\\
&=&((\pi _{pq^{-1}p^{-1}}\o \pi _{p^{-1}})\D (a))((i\o
 \pi _{p^{-1}})(\si (R)))R(1\o b).
 \end{eqnarray*}
\begin{flushright}
$\blacksquare$
\end{flushright}

{\bf 3.8 Lemma.} Take any $p, q\in G$. Let $a\in A_{p}$ and $b\in
 A_q$. Then we have

$$
\D (\widetilde{u})((i\o \pi _{q^{-1}})(\si (R)))((\pi _{p^{-1}}\o
\pi _{q^{-1}})(R))(a\o b)=(\widetilde{u}\o \widetilde{u})(a\o b).
$$

{\bf Proof.} We write $a\o b=\s \D (a_i)(1\o b_i)$ where
 $a_i\in A_{pq}$ and $b_i\in A_q$. By using Lemma 3.7, we become
 that

\begin{eqnarray*}
&&((i\o \pi _{p^{-1}})(\si (R)))R(a\o b)=\\
&&\quad \s ((\pi _{pq^{-1}p^{-1}}\o \pi _{p^{-1}})\D (a_i))
 ((i\o \pi _{p^{-1}})(\si (R)))R(1\o b_i).
\end{eqnarray*}

Therefore,

\begin{eqnarray*}
&& ((\pi _{pq^{-1}p^{-1}}\o \pi _{p^{-1}})\D (\widetilde{u}))
 ((i\o \pi _{p^{-1}})(\si (R)))R(a\o b)\\
&=& \s ((\pi _{pq^{-1}p^{-1}}\o \pi _{p^{-1}})\D
(\widetilde{u}a_i))
 ((i\o \pi _{p^{-1}})(\si (R)))R(1\o b_i)\\
&=& \s ((\pi _{pq^{-1}p^{-1}}\o \pi _{p^{-1}})\D (\pi _{pq}
(S(R^{(2)}))R^{(1)}a_i))((i\o \pi _{p^{-1}})(\si (R)))R(1\o b_i)\\
&=& \s ((\pi _{p}\o \pi _{q})(\D (S(R^{(2)}))))
\underline{((\pi_{pq^{-1}p^{-1}}\o \pi _{p^{-1}})(\D
(R^{(1)}a_i)))}\\
&&\quad \quad \quad \quad \quad \quad \quad \quad \quad
 \quad \quad \quad \quad \quad \quad \quad \quad
 \quad \quad \quad \quad \quad \quad
\underline{((i\o \pi _{p^{-1}})(\si (R)))R(1\o b_i)}
\end{eqnarray*}

By using Lemma 3.7 to the underlined expression, the above formula
 becomes

$$
\s ((\pi _{p}\o \pi _{q})(\D (S(R^{(2)}))))
 ((i\o \pi _{p^{-1}})(\si (R)))R\D (R^{(1)}a_i)(1\o b_i).
$$

 In the following, we use the letters $R, U, T, V, W,
 Z$ to denote the generalized $\pi $-matrix for $A$.
 The foregoing expression equals
\begin{eqnarray*}
&&\s ((\pi _{p}\o \pi _{q})(\D (S(R^{(2)}V^{(2)}))))
 (U^{(2)}\o \pi _{p^{-1}}(U^{(1)}))(T^{(1)}\o T^{(2)})
(\pi _q(R^{(1)})\o V^{(1)})\\
&&\quad \quad \quad \quad \quad \quad \quad \quad \quad \quad
\quad \quad \quad \quad \quad \quad \quad \quad \quad \quad \quad
\, \quad \quad \quad \quad \quad \quad \quad \quad \quad \quad
(a\o b)\\
&=& \s (\pi _{p}(S(W^{(2)}Z^{(2)}))\o \pi _{q}
 (S(R^{(2)}V^{(2)})))(U^{(2)}T^{(1)}\o \pi
 _{p^{-1}}(U^{(1)})T^{(2)})\\
 && \quad \quad \quad \quad \quad \quad \quad \quad
 \quad \quad \quad \quad \quad \quad \quad \quad
  \quad \quad \quad \quad
 (\pi _q(W^{(1)}R^{(1)}) \o Z^{(1)}V^{(1)})(a\o b)\\
&=&\s (S(\pi_{p}
(\underline{W^{(2)}Z^{(2)})})U^{(2)}\underline{T^{(1)}
 \pi _q(W^{(1)})}\pi _q(R^{(1)})a)\o \\
 &&\quad \quad \quad \quad \quad \quad \quad \quad
 \quad \quad \quad \quad \quad \quad \quad \quad \quad
 (\pi _{q} (S(R^{(2)}V^{(2)}))\pi _{p^{-1}}(U^{(1)})
 \underline{T^{(2)} Z^{(1)}}V^{(1)}b)\\
 &&\\
&&\mbox {We apply the $\pi $-Yang Baxter equation to the
 underlined expressions in the}\\
&& \mbox{sense of Proposition 3.3.}\\
&&\\
 &=& \s (S(\pi
 _{p}(W^{(2)}Z^{(2)}))U^{(2)}Z^{(1)}T^{(1)}
 \pi _q(R^{(1)})a)\o \\
 && \quad \quad \quad \quad \quad \quad \quad \quad
 \quad \quad \quad \quad \quad \quad \quad \quad
 (\pi _{q}
 (S(R^{(2)}V^{(2)}))\pi _{p^{-1}}(U^{(1)})W^{(1)}
 T^{(2)} V^{(1)}b).
\end{eqnarray*}

\quad \quad \quad  \quad \quad \quad  \quad \quad \quad \quad We
 notice that for all
 elements $x\in A$, we have
\begin{eqnarray*}
\quad \quad \quad \quad \quad \quad \quad \quad
&&\s \pi _{p^{-1}}(U^{(1)})W^{(1)}x\o S(\pi _p(W^{(2)}))U^{(2)}\\
&=&\s U^{(1)}W^{(1)}x\o S(\pi _p(W^{(2)}))\pi _p(U^{(2)})\\
&=&(i\o \pi _p)(\s U^{(1)}W^{(1)}x\o S(W^{(2)})U^{(2)})\\
&=&(i\o \pi _p)(i\o S)(\s U^{(1)}W^{(1)}x\o S^{-1}(U^{(2)})W^{(2)})\\
&=&(i\o \pi _p)(i\o S)(x\o 1)=x\o 1.
\end{eqnarray*}

Therefore, the above expression  can be simplified
 as follows

\begin{eqnarray*}
&&\s S(\pi _{p}(Z^{(2)}))Z^{(1)}T^{(1)}\pi _q(R^{(1)})a
 \o \pi _{q}(S(R^{(2)}V^{(2)}))T^{(2)}V^{(1)}b\\
&=& \s S(\pi _{p}(Z^{(2)}))Z^{(1)}a\o \pi _{q}
 (S(V^{(2)}))V^{(1)}b=(\widetilde{u}\o \widetilde{u})(a\o b).
\end{eqnarray*}

So we have proven that for all $a\in A_p$ and $b\in A_q$, we have

$$
((\pi _{pq^{-1}p^{-1}}\o \pi _{p^{-1}})\D (\widetilde{u}))((i\o
\pi _{p^{-1}})(\si (R)))R(a\o b)=
  (\widetilde{u}\o \widetilde{u})(a\o b).
$$

The above equation can be equivalently written as

$$
 \D (\widetilde{u})((\pi _{pqp^{-1}}\o i)(\si (R)))
((\pi _{pqp^{-1}}\o \pi _p)(R)) (\pi _{pqp^{-1}}(a)\o \pi _p(b))
=(\widetilde{u}\o \widetilde{u})(\pi _{pqp^{-1}}(a)\o \pi
 _p(b)).
$$

As $\pi $ is a crossing of $G$ on $A$, this latter
 equation can be written as follows. For all
 $a\in A_p$ and $b\in A_q$, we have

$$
 \D (\widetilde{u})((i\o \pi _{q^{-1}})(\si (R)))
((\pi _{p^{-1}}\o \pi _{q^{-1}})(R))(a\o b)= (\widetilde{u}\o
\widetilde{u})(a\o b).
$$

\begin{flushright}
$\blacksquare$
\end{flushright}

{\bf 3.9 Lemma.} Take any $p, q\in G$. Let $a\in A_p$ and $b\in
 A$. Then we have

$$
((\pi _p\o i)(\si (R)))R\D (\widetilde{u})(a\o b)
=(\widetilde{u}\o \widetilde{u})(a\o b).
$$
\\

{\bf Proof.} Following Lemma 3.7, we have for all
 $x\in A_{pq}$ and $b\in A_q$
\begin{eqnarray*}
&&((\pi _{pqp^{-1}}\o i)\si (R))((\pi _{pqp^{-1}}\o \pi _p)(R))
((\pi _{pqp^{-1}}\o \pi _p)(\D (x)))(1\o \pi _p(b))\\
&=&\D (x)((\pi _{pqp^{-1}}\o i)\si (R))((\pi _{pqp^{-1}}\o \pi
 _p)(R))(1\o \pi _p(b)).
\end{eqnarray*}

The above equation can also be written as follows. For all
 $x\in A_{qp}$ and $b\in A_q$, we have
\begin{eqnarray*}
&&((\pi _{q}\o i)\si (R))((\pi _{q}\o \pi _p)(R))
((\pi _{q}\o \pi _p)(\D (x)))(1\o b)\\
&=&\D (x)((\pi _{q}\o i)\si (R))((\pi _{q}\o \pi
 _p)(R))(1\o b).
\end{eqnarray*}

Therefore, we have for all $x\in M(A), b\in A_q$ and
  $a\in A_{qpq^{-1}}$
\begin{eqnarray*}
&&(\pi _q\o \pi _p)(((\pi _{p}\o i)(\si (R)))R\D
 (x))(a\o b)\\
&& \quad  =\D (x)((i\o \pi _{q^{-1}})(\si (R)))((\pi
 _{p^{-1}}\o \pi _{q^{-1}})(R))(a\o b).
\end{eqnarray*}

We set $x=\widetilde{u}$ and use Lemma 3.8. Then we have

$$
(\pi _q\o \pi _p)(((\pi _p\o i)(\si (R)))R\D (\widetilde{u})) (a\o
 b)=(\widetilde{u}\o \widetilde{u})(a\o b).
$$

Now, the formula in the statement easily follows. $\blacksquare$
\\
\\

{\bf 3.10 Proposition.} The multiplier
 $\widetilde{u}S(\widetilde{u})^{-1}$ is
 grouplike in $M(A)$.
\\

{\bf Proof.} First, we calculate the multiplier $\D
 (S(\widetilde{u}))$
 in $M(A\o A)$ by making use of Lemma 3.9.
 Let $p, q$ be in $G$ and take $a\in A_p$ and $b\in A_q$.

\begin{eqnarray*}
&&\D (S(\widetilde{u}))(a\o b)=((\si \circ (S\o S))\D
(\widetilde{u}))(a\o b)\\
&=&((\si \circ (S\o S))(R^{-1}((\pi _{q^{-1}}\o
 i)(\si (R^{-1})))(\widetilde{u}\o \widetilde{u})))(a\o b)\\
&=&(S(\widetilde{u})\o S(\widetilde{u}))((i\o \pi _{q^{-1}})
 ((S\o S)(R^{-1})))((S\o S)\si (R^{-1}))(a\o b).
\end{eqnarray*}

Therefore, we have

\begin{eqnarray*}
&&\D (S(\widetilde{u})^{-1})(a\o b)\\
&=&((S\o S)\si (R))((i\o \pi _{q^{-1}})
 ((S\o S)(R)))(S(\widetilde{u})^{-1}
 \o S(\widetilde{u})^{-1})(a\o b).
\end{eqnarray*}

As $(\widetilde{S}\o S)(R)=R$, we have for all $a\in A_p$
 and $b\in A_q$

\begin{eqnarray*}
&&\D (S(\widetilde{u})^{-1})(a\o b)\\
&=&((i\o \pi _{q^{-1}})(\si (R)))((\pi _{p^{-1}}\o \pi _{q^{-1}})
 (R)))(S(\widetilde{u})^{-1}\o S(\widetilde{u})^{-1})(a\o b).
\end{eqnarray*}

Take $p, q$ in $G$ and $a\in A_p, b\in A_q$, we combine
 the above result with Lemma 3.8 to obtain

\begin{eqnarray*}
\, \, \, \, \, \, \quad &&\D (\widetilde{u}S(\widetilde{u})^{-1})(a\o b)\\
&=&(\widetilde{u}\o \widetilde{u})((\pi _{p^{-1}}\o \pi
 _{q^{-1}})(R^{-1}))((i\o \pi _{q^{-1}})(\si (R^{-1})))
 ((i\o \pi _{q^{-1}})(\si (R)))\\
&&\quad \quad \quad \quad \quad
 \quad \quad \quad \quad \quad
((\pi _{p^{-1}}\o \pi _{q^{-1}})
 (R))(S(\widetilde{u})^{-1}\o S(\widetilde{u})^{-1})(a\o b)\\
&=&(\widetilde{u}S(\widetilde{u})^{-1}\o
\widetilde{u}S(\widetilde{u})^{-1})(a\o b).
\end{eqnarray*}
So we have proven that $\widetilde{u}S(\widetilde{u})^{-1}$
 is grouplike in $M(A)$. $\blacksquare$

\vskip 0.6cm

{\bf 3.11 The Drinfel'd double construction for $G$-cograded
 multiplier Hopf algebras}

\vskip 0.6cm

Let $\<A, B\>$ be a pairing of multiplier Hopf algebras. Suppose
 that $B$ is $G$-cograded. We put $B=\bigoplus _{p\in G}B_p$. Let
 $\pi $ denote a crossing of $G$ on $B$. This means that
 for all $p, q\in G$ we have $\pi _p(B_q)=B_{pqp^{-1}}$.

 The Drinfel'd double construction, denoted by $D^{\pi }$,
  is reviewed in the preliminaries. We have that $D^{\pi }$ is again a
 a $G$-cograded multiplier Hopf algebra with a nontrivial
 crossing of $G$ on $D^{\pi }$, see [De-VD3, Prop. 3.13].
 More precisely, we put $D^{\pi }_p=A\bowtie B_{p^{-1}}$
 for all $p\in G$. Define the automorphism $\pi '_p$
 on $A$ by the formula $\< \pi '_p(a), b\>=\< a, \pi _{p^{-1}}(b)\>$.
 The maps $\{ \pi '_p\o \pi _p, p\in G\}$ provide $D^{\pi }$ with a
 natural crossing of $G$ on $D^{\pi }$.

 In this section we investigate when $D^{\pi }$ is a
 $\pi $-quasitriangular $G$-cograded multiplier Hopf algebra in the sense
 of Definition 3.1.

  We suppose that the pairing $\< A, B\>$ has a canonical multiplier
  $W$ in $M(B\o A)$. This means that $W$ is invertible in $M(B\o A)$
  and for all $a\in A$ and $b\in B$, we have $\< W, a\o b\>=\< a, b\>$.
 From the definition, $W$ is unique in $M(B\o A)$. Following [De-VD2,
 Proposition 4.3], we have for all $a\in A$ and $b\in B$,
 $(\< a, \cdot \>\o i_A)(W)=a$ and $(i_B\o \< \cdot , b\>)(W)=b$.
  As $\< A_p, B_q\>=0$ whenever $p\ne q$,
  the following "covering" conditions on $W$ are quite natural.
  \\

  For all $b\in B_p$, we assume

$$
W(b\o 1)\in B_p\o A_p\quad \mbox{and}\quad (b\o 1)W\in B_p\o A_p.
$$

We will use the following notation
$$
W(b\o 1)=\s W^{(1)}b\o W^{(2)}\quad \mbox{and}\quad (b\o 1)W=\s
 bW^{(1)}\o W^{(2)}.
$$

Recall that $A$ and $B$ are embedded in $M(D^{\pi })$ in the
 following way

$$
A\lr M(D^{\pi }): \, a\mapsto a\bowtie 1 \quad \mbox {and} \quad
 B\lr M(D^{\pi }): \, b\mapsto 1\bowtie b.
$$

Therefore, we have the (non-degenerate) embedding

$$
 B\o A\lr M(D^{\pi }\o D^{\pi }): \,
 b\o a\mapsto (1\bowtie b)\o (a\bowtie 1).
$$

By extending this latter algebra embedding to $M(B\o A)$, it
 makes sense to consider $W\in M(B\o A)$ as a multiplier
 in $M(D^{\pi }\o D^{\pi })$. The following theorem generalizes
 [Zun, Theorem 5.6]
\\
\\

{\bf 3.11.1. Theorem.} Take the notations and assumptions as
 above. The Drinfel'd double $D^{\pi }$ of the pair $\< A, B\>$ is
 $\pi $-quasitriangular for the natural crossing of $G$ on $D^{\pi
 }$. A $\pi $-matrix is given by the embedding of the canonical
 multiplier of the pair $\< A, B\>$.

{\bf Proof.} (1) For all $p\in G$, we have
 $(\pi _p\o \pi _p')(W)=W$ in $M(B\o A)$. Therefore,
 we have

$$
((\pi '_p\o \pi _p)\o (\pi '_p\o \pi _p))(W)=W \quad \mbox
{in}\quad M(D^{\pi }\o D^{\pi }).
$$
\\

(2) We prove
$$
(\widetilde{\D} _D\o i_D)(W)=W^{13}W^{23}\quad \mbox{in}\, \,
 M(D^{\pi }\o D^{\pi }\o D^{\pi }).
$$

Recall that for all $a, x\in A, b\in B$ and $y\in B_q$, we have

$$
\widetilde{\D }(a\bowtie b)((1\bowtie 1)\o (x\bowtie y)) =\s (\pi
 '_q(a_{(2)})\bowtie b_{(1)})\o (a_{(1)}\bowtie b_{(2)})(x\bowtie y)
$$

Therefore,
$$
(\widetilde{\D} _D\o i_D)(W)=(\D _B\o i_A)(W)\quad \mbox{in}\,
 \, \, M(D^{\pi }\o D^{\pi }\o D^{\pi })
$$

By [De-VD2, Proposition 4.4], this expression equals
 $W^{13}W^{23}$.
\\

(3) We prove

$$
(i_D\o \D _D)(W)=W^{13}W^{12}\quad \mbox{in}\, \,
 M(D^{\pi }\o D^{\pi }\o D^{\pi }).
$$

Indeed,

$$
(i_D\o \D _D)(W)=(i_B\o \D ^{cop}_A)(W)\quad \mbox{in}\,
 \, \, M(D^{\pi }\o D^{\pi }\o D^{\pi })
$$

Again by [De-VD2, Proposition 4.4], this expression equals
 $W^{13}W^{12}$.
\\

(4) We prove that for all $a\in A$ and $b\in B$

$$
W\D (a\bowtie b)=\widetilde{\D }^{cop}(a\bowtie b)W\quad \,
 \mbox {in}\quad M(D^{\pi }\o D^{\pi }).
$$

Take $x, a\in A, b\in B_{pq}, y\in B_p$ and $y'\in B_q$.

\begin{eqnarray*}
&& W\D (a\bowtie b)((1\bowtie y')\o (x\bowtie y))\\
&=&\s ((1\bowtie W^{(1)})(a_{(2)}\bowtie \pi
 _{p^{-1}}(b_{(1)})y'))\o
 ((W^{(2)}a_{(1)}\bowtie b_{(2)})(x\bowtie y)).
\end{eqnarray*}

Again use [De-VD2, Proposition 4.4] and let the canonical
 multiplier be denoted by the letters $W, T, V$.
 Then the expression above can be written as

\begin{eqnarray*}
&&\s (\< a_{(2)}, S^{-1}(T^{(1)})\>\<
 a_{(4)}, \pi _{q^{-1}}(W^{(1)})\>(a_{(3)}\bowtie V^{(1)}\pi
 _{p^{-1}}(b_{(1)})y') \o \\
 &&\quad \quad \quad \quad \quad
 \quad \quad \quad \quad \quad
 \quad \quad \quad \quad \quad
 \quad \quad \quad \quad \quad
 ((W^{(2)}V^{(2)}T^{(2)}a_{(1)}\bowtie b_{(2)})(x\bowtie y))\\
&&=(a_{(1)}\bowtie V^{(1)}\pi _{p^{-1}}(b_{(1)})y')\o
 ((\pi '_q(a_{(2)})V^{(2)}\bowtie b_{(2)})(x\bowtie y)).
\end{eqnarray*}
\\
\\

Again, take $a, x\in A, b\in B_{pq}, y\in B_p$ and $y'\in B_q$.
 We calculate the expression
\begin{eqnarray*}
&&\widetilde{\D }^{cop}(a\bowtie b)W((1\bowtie y')\o (x\bowtie
 y))\\
&=&\widetilde{\D }^{cop}(a\bowtie b)((1\bowtie W^{(1)}y')\o
 (W^{(2)}x\bowtie y))\\
&=&\s (a_{(1)}\bowtie b_{(2)}W^{(1)}y')\o
 ((\pi '_q(a_{(2)})\bowtie  b_{(1)})(W^{(2)}\bowtie 1)
 (x\bowtie y))\\
&=&\s \<W^{(2)}, S^{-1}(b_{(3)})\>\< T^{(2)}, \pi
 _{p^{-1}}(b_{(1)})\>(a_{(1)}\bowtie
 b_{(4)}W^{(1)}V^{(1)}T^{(1)}y')\o \\
 &&\quad \quad \quad \quad \quad
 \quad \quad \quad \quad \quad
 \quad \quad \quad \quad \quad
 \quad \quad \quad \quad \quad
 \quad
 ((\pi '_q(a_{(2)})V^{(2)}\bowtie  b_{(2)})(x\bowtie y))\\
&=&\s (a_{(1)}\bowtie V^{(1)}\pi
 _{p^{-1}}(b_{(1)})y')\o
 ((\pi '_q(a_{(2)})V^{(2)}\bowtie  b_{(2)})(x\bowtie y)).
  \end{eqnarray*}

So, we have proven for all $a\in A$ and $b\in B$

$$
W\D (a\bowtie b)=\widetilde{\D }^{cop}(a\bowtie b)W\quad
 \mbox{in}\, \, M(D^{\pi }\o D^{\pi }).
$$

We now conclude that the embedding of $W$ in $M(D^{\pi }\o
 D^{\pi })$ makes $D^{\pi }$ into a $\pi $-quasitriangular
 multiplier Hopf algebra, in the sense of Definition 3.1.
 $\blacksquare$
\\
\\

{\bf 3.11.2.  "Finite type" Hopf group-coalgebras} \quad
 We give concrete data to illustrate Theorem 3.11.1.
 Let $G$ be any  group. Consider a "finite
 type" Hopf $G$-coalgebra as given in [T-Section 11].
  In [A-De-VD, Theorem 1.5]
 we have shown how to consider any  Hopf group-coalgebra into
 the framework of $G$-cograded multiplier Hopf algebras.
 In this framework, a "finite type" Hopf $G$-coalgebra is
 dealt as a $G$-cograded multiplier
 Hopf algebra $A=\bigoplus _{p\in G}A_p$ where each
 $A_p$ is a finite-dimensional algebra with a unit.
 A crossing in the sense of [T-Section 11] defines a crossing
 on $A$. In [A-De-VD, Theorem 3.9] we have constructed
 an integral on $A$. Therefore, we can consider the
  reduced dual multiplier Hopf
 algebra which is given by the Hopf algebra
 $A^*=\bigoplus _{p\in G}(A_p)'$, see [A-De-VD, Corollary 3.4].

Consider the pair $\< A^*, A\>$. Let $\{ f_{pi}\}\subset (A_p)'$
 and $\{ e_{pi}\}\subset A_p$ be dual bases. Then the canonical
 multiplier of the pair $\< A^*, A\>$ is given as

 $$
W=\s _{p,i}e_{pi}\o f_{pi} \quad \mbox {in} \, \, M(A\o A^*).
 $$

Clearly, $W(a\o 1)$ and $(a\o 1)W$ are in $A\o A^*$ for all
 $a\in A$.  By Theorem 3.11.1, we obtain that the
 embedding of $W$ in $M(D^{\pi }\o
 D^{\pi })$ is a $\pi $-quasitriangular structure for $D^{\pi }$,
 considered for the natural crossing $\{ \pi '_p\o \pi _p, p\in
 G\}$. The embedding of $W$ in $M(D^{\pi }\o D^{\pi })$ is
 considered as $W=\s _{p\in G}(1\bowtie e_{pi})\o (f_{pi}\bowtie
 1)$. This result is also in [Zun, Theorem 5.6].

Take the "finite type" Hopf group-coalgebras where
 $dim (A_p)=1$ for all $p\in G$.
 The $G$-cograded multiplier Hopf algebra $A$ is given as
 $A=K(G)$, the multiplier Hopf algebra of
  all complex-valued functions with finite support in $G$,
  see Preliminaries. The reduced dual multiplier Hopf
  algebra $A^*$ is now given by the usual group Hopf algebra
  $A^*=k[G]$. The canonical
  multiplier of the pair $\< k[G], K(G)\>$
  is given as $W=\s _{p\in G}\d _p\o
  u_p$ in $M(K(G)\o k[G])$.
  The crossing of $G$ on $K(G)$ is given as
  $\pi _p(\d _q)=\d _{pqp^{-1}}$ for all $p, q\in G$.
 We have $\widetilde{K(G)}=K(G)^{cop}$. The Drinfel'd
 double is just a usual tensor product
 $D^{\pi }=k[G]\o K(G)^{cop}$.
 An easy calculation shows that
 $\widetilde{u}=\s _p
 u_{p^{-1}}\o \d _p$ in  $M(D^{\pi })$.
 We have in this case
 $\widetilde{u}S(\widetilde{u})^{-1}=u_e\o \s _{p\in G}\d _p=1$
  in $M(D^{\pi })$.
\\

\vskip 1.5cm

\begin{center}
 {\bf ACKNOWLEDGEMENT}
\end{center}

  The third author was supported by the Research Council of the
   K.U.Leuven. He is very grateful to the research group of
   K.U.Leuven for providing a good atmosphere to work.

\vskip 1.5cm

\begin{center}
{\bf REFERENCES}
\end{center}

[A-De-VD] A.T.Abd El-hafez, L. Delvaux, and A. Van Daele,
  Group-cograded multiplier Hopf ($*$-)algebra,
   math. QA/0404026.
\smallskip

[B-D-G-N] M. Beattie, S. D$\check{a}$sc$\check{a}$lescu, L.
 Gr$\ddot{u}$nenfelder, C. N$\check{a}$st$\check{a}$sescu,
 Finiteness conditions, Co-Frobenius Hopf algebras
 and Quantum groups, J. Algebra 200(1998), 312-333.
\smallskip

[De] L. Delvaux, Pairing and Drinfel'd double of Ore-extensions,
 Comm. in Algebra 29(7)(2001), 3167-3177.
\smallskip

[De-VD1] L. Delvaux and A. Van Daele, The Drinfel'd double
 of multiplier Hopf algebras, J. Algebra 272(1)(2004), 273-291.
\smallskip

[De-VD2] L. Delvaux and A. Van Daele, The Drinfel'd double versus
 the Heisenberg double for an algebraic quantum group,
 J. Pure and Appl. Algbera, 190 (2004), 59-84.
\smallskip

[De-VD3] L. Delvaux and A. Van Daele, The Drinfled double for
 group-cograded multiplier Hopf algebra. Preprint LUC and K.U. Leuven
 \smallskip

[Dra-VD] B. Drabant and A. Van Daele, Pairing and the quantum
 double of multiplier Hopf algebras, Algebras and Representation Theory
 4 (2001), 109-132.
\smallskip

[Dri1] V. G. Drinfel'd, Quantum groups, Proc. Int. Cong. Math.,
 dr 86 Berkeley, 1(1986), 789-820.
\smallskip

[Dri2] V. G. Drinfel'd, On almost cocommutative Hopf algebras,
 Leningrad Math. J. 1(1990), 321-342.
\smallskip

[Ra] D. E. Radford, On the antipode of a quasitriangular Hopf
 algebra, J. Algebra 151 (1992), 1-11.
\smallskip

[T] V. G. Turaev,  Homotopy field theory in dimension $3$ and
 crossed group-categories. Preprint GT/0005291.
\smallskip

[VD-1] A. Van Daele, Multiplier Hopf algebras, Trans. Soc.
 342(2)(1994), 917-932.
\smallskip

[VD-2] A. Van Daele, An algebraic framework for group duality,
 Advances in Mathematics, 140(1998), 323-366.

[VD-Z] A. Van Daele and Y. Zhang, Multiplier Hopf algebras of
 discrete type, J. Algebra 214(1999), 400-417.
\smallskip

 [Vir] A. Virelizier, Hopf group-coalgebras, J. Pure and Applied
 Algebra, (2002)171, 75-122.
\smallskip

[Z] Y. Zhang, The quantum Double of a co-Frobenius Hopf algebra,
 Comm. in Algebra 27(3) (1999), 1413-1427.
\smallskip

[Zun] M., Zunino, Double construction for crossed Hopf
 coalgebra. J. Algebra 278(2004), 43-75.

 \end{document}